\magnification=1200

\loadmsam
\loadmsbm
\loadeufm
\loadeusm
\UseAMSsymbols

\font\BIGtitle=cmr10 scaled\magstep3
\font\bigtitle=cmr10 scaled\magstep1
\font\boldsectionfont=cmb10 scaled\magstep1
\font\section=cmsy10 scaled\magstep1

\def\scr#1{{\fam\eusmfam\relax#1}}

\def\scrB{{\scr B}}

\def\scrL{{\scr L}}

\def\scrO{{\scr O}}

\def\scrW{{\scr W}}

\def\gr#1{{\fam\eufmfam\relax#1}}

	\def\grg{{\gr g}}

	\def\grn{{\gr n}}

	\def\grs{{\gr s}}

\def\db#1{{\fam\msbfam\relax#1}}

 \def\dbF{{\db F}}
\def\dbG{{\db G}}

 \def\dbN{{\db N}}

\def\eps{{\varepsilon}}

\def\im{\text{Im}}

\def\Ker{\text{Ker}}
\def\der{\text{der}}

\def\sh{\hbox{\rm sh}}

\def\ab{\text{ab}}

\def\ad{\text{ad}}

\def\sh{\text{sh}}

\def\Spec{\text{Spec}}

\def\red{\text{red}}

\def\n{\text{n}}

\def\sc{\text{sc}}

\def\Lie{\text{Lie}}

\def\leaderfill{\leaders\hbox to 1em
     {\hss.\hss}\hfill}
\def\nspace{\lineskip=1pt\baselineskip=12pt\lineskiplimit=0pt}

\def\finishproclaim{\par\rm
     \ifdim\lastskip<\medskipamount\removelastskip
     \penalty55\medskip\fi}
\def\endproof{$\hfill \square$}
\def\proof{\par\noindent {\it Proof:}\enspace}
\def\references#1{\par
  \centerline{\boldsectionfont References}\medskip
     \parindent=#1pt\nspace}
\def\Ref[#1]{\par\hang\indent\llap{\hbox to\parindent
     {[#1]\hfil\enspace}}\ignorespaces}
\def\Item#1{\par\smallskip\hang\indent\llap{\hbox to\parindent
     {#1\hfill$\,\,$}}\ignorespaces}
\def\ItemItem#1{\par\indent\hangindent2\parindent
     \hbox to \parindent{#1\hfill\enspace}\ignorespaces}

\def\Le{{\mathchoice{\,{\scriptstyle\le}\,}
  {\,{\scriptstyle\le}\,}
  {\,{\scriptscriptstyle\le}\,}{\,{\scriptscriptstyle\le}\,}}}
\def\Ge{{\mathchoice{\,{\scriptstyle\ge}\,}
  {\,{\scriptstyle\ge}\,}
  {\,{\scriptscriptstyle\ge}\,}{\,{\scriptscriptstyle\ge}\,}}}

\def\arrowsim{\,\smash{\mathop{\to}\limits^{\lower1.5pt
  \hbox{$\scriptstyle\sim$}}}\,}

\def\doublemaprights#1#2#3#4{\raise3pt\hbox{$\mathop{\,\,\hbox to     
#1pt{\rightarrowfill}\kern-30pt\lower3.95pt\hbox to
     #2pt{\rightarrowfill}\,\,}\limits_{#3}^{#4}$}}

\def\rightcapdownarrow{\raise9pt\hbox{$\ssize\cap$}\kern-7.75pt
     \Big\downarrow}

\def\rcapmapdown#1{\rightcapdownarrow\kern-1.0pt\vcenter{
     \hbox{$\scriptstyle#1$}}}

\def\rmapdown#1{\Big\downarrow\kern-1.0pt\vcenter{
     \hbox{$\scriptstyle#1$}}}
\def\rightsubsetarrow#1{{\ssize\subset}\kern-4.5pt\lower2.85pt
     \hbox to #1pt{\rightarrowfill}}
\def\longtwoheadedrightarrow#1{\raise2.2pt\hbox to #1pt{\hrulefill}
     \!\!\!\twoheadrightarrow}

\NoBlackBoxes
\parindent=25pt
\document
\footline={\hfil}

\null
\noindent 
\centerline{\BIGtitle On two theorems for flat, affine group schemes}
\centerline{\BIGtitle over a discrete valuation ring}

\vskip 0.3 cm
\centerline{\bigtitle Adrian Vasiu, UA, July 27, 2004}
\footline={\hfill}
\vskip 0.3 cm
\noindent
{\bf ABSTRACT.} We include short and elementary proofs of two theorems that characterize reductive group schemes over a discrete valuation ring, in a slightly more general context. 


$\backslash$


\footline={\hss\tenrm \folio\hss}
\pageno=1

\bigskip
\noindent
{\boldsectionfont \S1. Introduction}
\bigskip

Let $k$ be a field. Let $p\in\{0\}\cup\{n\in\dbN|n\;\text{is}\;\text{a}\;\text{prime}\}$ be the characteristic of $k$. Let $V$ be a discrete valuation ring of residue field $k$. Let $\pi$ be a uniformizer of $V$ and let $K:=V[{1\over\pi}]$ be the field of fractions of $V$. Let $F=\Spec(P)$ and $G=\Spec(R)$ be flat, affine group schemes over $V$. We will assume that $F$ is a reductive group scheme over $V$; so $F$ is smooth over $V$ and its fibres are connected and reductive groups over fields. In this paper we present elementary and short proofs of the following two basic theorems on reductive group schemes over $V$. 

\medskip\smallskip\noindent
{\bf 1.1. Theorem.} {\it Let $f:F\to G$ be a homomorphism such that its generic fibre $f_K:F_K\to G_K$ is a closed embedding. Then the following four properties hold:

\medskip
{\bf (a)} Suppose $f_K$ is an isomorphism and $G$ is of finite type. Then $F$ is the smoothening of $G$ in the sense of [BLR, Thm. 5, p. 175] (the definition is recalled in 2.3.2).

\smallskip
{\bf (b)} The kernel $\Ker(f_k:F_k\to G_k)$ is a unipotent, connected group of dimension $0$.

\smallskip
{\bf (c)} The homomorphism $f$ is finite. 

\smallskip
{\bf (d)} If $p=2$, then we assume that $F_{\bar K}$  has no normal subgroup that is an $SO_{2n+1}$ group for some $n\in\dbN$. Then $f$ itself is a closed embedding.}

\medskip\smallskip\noindent
{\bf 1.2. Theorem.} {\it Suppose $G_K$ is a connected, smooth group over $K$ and the identity component $(G_k)_{\red}^0$ of the reduced group $(G_k)_{\red}$ is a reductive group over $k$ of the same dimension as $G_K$. Then the following three properties hold:

\medskip
{\bf (a)} The group scheme $G$ is of finite type over $V$, has a connected special fibre $G_k$, and its generic fibre $G_K$ is a reductive group. 

\smallskip
{\bf (b)} If $p=2$, we assume that the group $G_{\bar K}$ has no normal subgroup that is isomorphic to $SO_{2n+1}$ for some $n\in\dbN$. Then $G$ is a reductive group scheme over $V$.

\smallskip
{\bf (c)} The normalization $G^{\n}$ of $G$ is a finite $G$-scheme. Moreover, there is a faithfully flat $V$-algebra $\tilde V$ that is a discrete valuation ring and such that the normalization $(G_{\tilde V})^{\n}$ of $G_{\tilde V}$ is a reductive group scheme over $\tilde V$.}

\medskip
If $V$ is of mixed characteristic $(0,p)$, then 1.1 (a) and (b) were first proved in [Va1, proof of 3.1.2.1 c)]. For the sake of completeness, in \S3 we recall the proofs of 1.1 (a) and (b) in a slightly enlarged manner that recalls elementary material. The passage from 1.1 (a) and (b) to 1.1 (d) (resp. 1.1 (c)) is a direct consequence of the classification of the ideals of the Lie algebras of adjoint and simply connected semisimple groups over $\bar k$ (resp. is only a variant of Zariski Main Theorem). In [Va1, 3.1.2.1 c)] we overlooked the phenomenon of exceptional nilpotent such ideals and so the extra hypothesis of 1.1 (d) for $p=2$ does not show up in [Va1, 3.1.2.1 c)]. Theorem 1.1 (resp. Theorem 1.2) is also proved in [PY], under the extra assumption that $G_K$ is of finite type over $K$ (resp. that $G_k$ is of finite type over $k$). The methods of [PY] are essentially the same as of [Va1] except that the  proofs of [PY] contain many unneeded parts. Here we follow [Va1] to get short and efficient proofs of 1.1 and 1.2. The importance of 1.1 (resp. 1.2) stems from its fundamental applications to integral canonical models of Shimura varieties of Hodge (resp. to dual groups), cf. [Va1] (resp. cf. [PY]). We also mention that 1.1 (c) and (d) are powerful tools in extending results on semisimple groups over a field, from characteristic $0$ to arbitrary positive characteristic. As an exemplification, in \S5 we include such an application that pertains to adjoint groups.  

The proofs of 1.1 and 1.2 are carried on in \S3 and \S4 (respectively). Few notations and preliminaries needed in \S3 and \S4 are gathered in \S2. We would like to thank U of Arizona for good conditions for the writing of this work. We would like to thank G. Prasad for pointing out that [Va1, 3.1.2.1 c)] omits to add the extra hypothesis of 1.1 (d) for $p=2$.

\bigskip
\noindent
{\boldsectionfont \S2. Preliminaries}
\bigskip 
In 2.1 we list notations. Elementary properties of Lie algebras, dilatations, representation theory, and quasi-sections are recalled in 2.2 to 2.5. In 2.6 we present in an accessible way a result of the classical Bruhat--Tits theory for reductive groups over $K$.

\medskip\smallskip\noindent
{\bf 2.1. Notations.} 
Let $k$, $V$, $K$, $\pi$, $F=\Spec(P)$, and $G=\Spec(R)$ be as in \S1. Let $V^{\sh}$ be the strict henselization of $V$. Let $\scrW$ be the set of finite, discrete valuation ring extensions of  the completion of $V^{\sh}$. If $H$ is a reductive group scheme over an affine scheme $\Spec(A)$, let $Z(H)$, $H^{\der}$, $H^{\ad}$, and  $H^{\ab}$ denote the center, the derived group, the adjoint group, and respectively the abelianization of $H$. So $Z^{\ab}(H)=H/H^{\der}$ and $H^{\ad}=H/Z(H)$. Let $H^{\sc}$ be the simply connected semisimple group cover of $H^{\der}$. For a free $A$-module $N$ of finite rank, let $GL(N)$ be the group scheme over $\Spec(A)$ of linear automorphisms of $N$. 

\medskip\smallskip\noindent
{\bf 2.2. Lie algebras.} Let  $\Lie(H)$ be the Lie algebra of an affine group $H$ of finite type over $k$ or over $V$. We view $\Lie(H)$ as a (left) $H$-module via the adjoint representation. Until 2.3 we assume that $k=\bar k$ and that $H$ is over $k$. For $x\in\Lie(H)$ we have a unique and functorial Jordan decomposition $x=x_s+x_n$ such that for any monomorphism $H\hookrightarrow GL_m$ ($m\in\dbN$), $x_s\in\Lie(GL_m)$ is semisimple, $x_n\in\Lie(GL_m)$ is nilpotent, and $x_s$ and $x_n$ are polynomials in $x\in \Lie(GL_m)$ (cf. [Bo, Ch. 1, \S 4]). We say $x$ is nilpotent (resp. semisimple) if $x=x_n$ (resp. $x=x_s$). So if $H$ is $\dbG_a$ (resp. $\dbG_m$), then $x$ is nilpotent (resp. semisimple). 

\medskip\noindent
{\bf 2.2.1. Lemma.} {\it Let $H$ be a reductive group over $k=\bar k$. Let $\grn$ be a non-zero ideal of $\Lie(H)$ that is formed by nilpotent elements and is a simple $H$-module. Then $p=2$ and $H$ has a normal subgroup $H_0$ isomorphic to $SO_{2n+1}$ ($n\in\dbN$) and such that $\grn\subset\Lie(H_0)$.}

\medskip
\proof
The image of $\grn$ in $\Lie(H^{\ab})$ is trivial. So $\grn\subset\Lie(H^{\der})$. So we can assume that $H=H^{\der}$ is semisimple. Let $\grn^{\ad}:=\im(\grn\to\Lie(H^{\ad}))$. As $\Lie(Z(H))$ is formed by semisimple elements, the Lie homomorphism $\grn\to \grn^{\ad}$ is an isomorphism. Thus $\grn^{\ad}$ is a non-zero ideal of $\Lie(H^{\ad})$ that is a simple $H^{\ad}$-module. Let $H^{\ad}=\prod_{i\in I} H_i$ be the product decomposition into simple groups. As $H^{\ad}$ is adjoint, there is no element of $\Lie(H^{\ad})$ fixed by $H^{\ad}$. So as $\grn$ is a simple $H^{\ad}$-module, we get that there is $i_0\in I$ such that $\grn^{\ad}\subset\Lie(H_{i_0})$. Thus we can assume $H^{\ad}$ is a simple group.  Let $T$ be a maximal torus of $H^{\ad}$ and let $\scrL$ be the Lie type of $H^{\ad}$. As $\Lie(H^{\ad})$ has non-zero semisimple elements, we have $\grn^{\ad}\neq\Lie(H^{\ad})$. So as $\grn^{\ad}$ is a simple $H^{\ad}$-module, from [Pi, Prop. 1.11] we get that either $\grn^{\ad}=\im(\Lie(H^{\sc})\to\Lie(H^{\ad}))$ or $(p,\scrL)\in\{(2,F_4),(3,G_2),(2,B_n),(2,C_n)|n\in\dbN\}$.

If $\grn^{\ad}=\im(\Lie(H^{\sc})\to\Lie(H^{\ad}))$, then $\grn^{\ad}$ has non-zero semisimple elements except when $\dim_k(\Lie(Z(H^{\sc})))$ is the rank of $\scrL$, i.e. except when $(p,\scrL)=(2,A_1)$ (cf. loc. cit.).

If $(p,\scrL)$ is $(2,F_4)$ (resp. $(3,G_2)$), then $\grn^{\ad}$ is the unique proper ideal of $\Lie(H^{\ad})$ and so it is generated by the direct sum $\grs$ of the eigenspaces of the adjoint action of $T$ on $\Lie(H^{\ad})$ that correspond to short roots. Thus $\dim_k(\grn^{\ad})$ is $26$ (resp. $7$) and $\grn^{\ad}\subset\grs\oplus\Lie(T)$, cf. [Hi, pp. 408--409] applied to a semisimple group over $\dbF_p$ whose extension to $k$ is $H^{\sc}$. But $\scrL$ has $24$ (resp. $6$) short roots, cf. [Bou, PLATES VIII and IX]. So $\dim_k(\Lie(T)\cap\grn^{\ad})$ is $2=26-24$ (resp. $1=7-6$) and so $\grn^{\ad}$ has non-zero semisimple elements. 
If $(p,\scrL)=(2,C_n)$ with $n\Ge 3$, then we similarly argue that $\grn^{\ad}$ is the ideal associated to short roots and that $\dim_k(\Lie(T)\cap\grn^{\ad})=2n^2-n-1-\eps-2(n^2-n)=n-1-\eps>1$, where $\eps$ is $1$ if $n$ is even and is $0$ if $n$ is odd (cf. [Hi, p. 409] and [Bou, PLATE III]). 

Thus we have $(p,\scrL)=(2,B_n)$, with $n\in\dbN$. As $\grn^{\ad}$ is a simple $H^{\ad}$-submodule of $\Lie(H^{\ad})$, from [Pi, Prop. 1.11] we get that that $\grn^{\ad}$ is the direct sum of the eigenspaces of the adjoint action of $T$ on $\Lie(H^{\ad})$ that correspond to short roots, that $\dim_k(\grn^{\ad})=2n$, and that $\grn^{\ad}\subset \im(\Lie(H^{\sc})\to\Lie(H^{\ad}))$. Let $\grn^{\sc}$ be the inverse image of $\grn$ in $\Lie(H^{\sc})$. We have $\dim_k(\grn^{\sc})=\dim_k(\grn)+\dim_k(\Lie(\Ker(H^{\sc}\to H)))=2n+\dim_k(\Lie(\Ker(H^{\sc}\to H)))$. But as $\grn^{\sc}$ is an $H^{\sc}$-module, we have $\Lie(Z(H^{\sc}))\subset\grn^{\sc}$ (cf. [Pi, Prop. 1.11]). Thus $\dim_k(\grn^{\sc})=2n+1$. Thus $\dim_k(\Lie(\Ker(H^{\sc}\to H)))=2n+1-2n=1$ and so $\Ker(H^{\sc}\to H)=Z(H^{\sc})\tilde\to\mu_{2}$. Thus $H$ is an adjoint group and so isomorphic to $SO_{2n+1}$.\endproof

\medskip\smallskip\noindent
{\bf 2.3. Dilatations.} In this section we assume $G_K$ is reduced. Let $S$ be a reduced subgroup of $G_k$. Let $J$ be the ideal of $R$ that defines $S$ and let $I_R$ be the ideal of $R$ that defines the identity section of $G$. Let $R_1$ be the $R$-subalgebra of $R[{1\over {\pi}}]$ generated by ${x\over\pi}$, where $x\in J$. By the dilatation of $G$ centered on $S$ one means the affine scheme $G_1:=\Spec(R_1)$; it has a canonical structure of a flat, affine group scheme over $V$ and the morphism $G_1\to G$ is a homomorphism whose special fibre factors through the closed embedding $S\hookrightarrow G_k$ (cf. [BLR, Prop. 1 and 2, pp. 63--64]). In \S3 we will need the following elementary Lemma.

\medskip\noindent
{\bf 2.3.1. Lemma.} {\it Suppose $G$ is a closed subgroup of a smooth group scheme $H$ over $V$ of relative dimension $l$. Then $\Ker(G_{1k}\to G_k)$ is isomorphic to a subgroup of $\dbG_a^{l-\dim(S)}$. Moreover, if $G=H$, then $\Ker(G_{1k}\to G_k)\arrowsim\dbG_a^{l-\dim(S)}$.}

\medskip
\proof
We can assume $V$ is complete. As $G_1$ is a closed subgroup of the dilatation of $H$ centered on $S$ (cf. [BLR, Prop. 2 (c) and (d), p. 64]), we can also assume that $G=H$. Let $\hat R$ and $\hat I_R$ be the completions of $R$ and $I_R$ with respect to the $I_R$-topology. Let $s:=\dim(S)$. As $G$ is smooth over $V$, we can write $\hat R=V\oplus \hat I_R=V[[x_1,...,x_{l-s},y_1,...,y_s]]$, where $x_1$, ..., $x_{l-s}$, $y_1$, ..., $y_s\in \hat I_R$  are such that the ideal $(x_1,...,x_{l-s},\pi)$ of $\hat R$ defines the completion of $S$ along its identity section. We have $R_1\otimes_{R} \hat R=V[[x_1,...,x_{l-s},y_1,...,y_s]][{x_1\over\pi},...,{x_{l-s}\over\pi}]$. Let $\delta:\hat R\to \widehat{\hat R\otimes_V  \hat R}$ be the comultiplication map of the formal Lie group of $G$.
As $S$ is a subgroup of $G_k$, for $i\in\{1,...,l-s\}$ we have $\delta(x_i)=x_i\otimes 1+1\otimes x_i+\sum_{j\in I_i^a} (a_{ij}\otimes b_{ij})+\sum_{j\in I_i^b} (a_{ij}\otimes b_{ij})$, where $I_i^a$ and $I_i^b$ are finite sets, where each $a_{ij}$ and $b_{ij}$ belong to $\hat I_R$,  and where for each $j\in I_i^a$ (resp. $j\in I_i^b$) the element $a_{ij}\in\hat I_R$ (resp. $b_{ij}\in\hat I_R$) is divisible by either some $x_u$ or by some $\pi y_v$; here $u\in\{1,...,l-s\}$ and $v\in\{1,...,s\}$.

We have $\Ker(G_{1k}\to G_k)=\Spec(A_1)$, where $A_1:=R_1\otimes_R \hat R/(x_1,...,x_{l-s},y_1,...,y_s)$. Let $\bar x_i$ be the image  of ${x_i\over{\pi}}$ in $A_1$. So $A_1$ is a $k$-algebra generated by $\bar x_1$, ..., $\bar x_{l-s}$. Taking the identity $\delta({x_i\over\pi})= {x_i\over\pi}\otimes 1+1\otimes {x_i\over\pi}+\sum_{j\in I_i^a} {a_{ij}\over \pi}\otimes b_{ij}+\sum_{j\in I_i^b} a_{ij}\otimes {b_{ij}\over\pi}$ modulo the ideal $(x_1,...,x_{l-s},y_1,...,y_s)$ of $R_1\otimes_R \hat R$, we get that the comultiplication map $\delta_1:A_1\to A_1\otimes_k A_1$ of the group $\Ker(G_{1k}\to G_k)$ is such that $\delta_1(\bar x_i)=\bar x_i\otimes 1+1\otimes\bar x_i$ for any $i\in\{1,...,l-s\}$. As $G_1$ is smooth over $V$ (cf. [BLR, Prop. 3, p. 64]) of relative dimension $l$, $\Ker(G_{1k}\to G_k)=\Ker(G_{1k}\to S)$ has dimension at least $l-s$. So as $A_1$ is $k$-generated by $\bar x_1$, ..., $\bar x_{l-s}$ and its dimension is at least $l-s$, we get that $A_1=k[\bar x_1,...,\bar x_{l-s}]$ is a polynomial $k$-algebra. From the description of $\delta_1$ we get that $\Ker(G_{1k}\to G_k)$ is isomorphic to $\dbG_a^{l-s}$.\endproof

\medskip\noindent
{\bf 2.3.2. Smoothening.} We assume that $G_K$ is smooth over $K$ and that $G$ is of finite type over $V$. We take $S$ to be the Zariski closure in $G_k$ of the special fibres of all morphisms $\Spec(V^{\sh})\to G$ of $V$-schemes. We refer to $G_1\to G$ as the canonical dilatation of $G$; it is a morphism of finite type. There is $m\in\dbN$ and a finite sequence of canonical dilatations $G^\prime:=G_m\to G_{m-1}\to ...\to G_1\to G_0:=G$ such that $G^\prime$ is uniquely determined by the following two properties (cf. [BLR, pp. 174--175]):

\medskip
{\bf (i)} the affine group scheme $G^\prime$ is smooth and of finite type over $V$; 

\smallskip
{\bf (ii)} if $Y$ is a smooth $V$-scheme and if $Y\to G$ is a morphism of $V$-schemes, then $Y\to G$ factors uniquely through the homomorphism $G^\prime\to G$. 

\smallskip
From very definitions, we get that $G^\prime(V)=G(V)$ and that the smoothening $(G_{V^{\sh}})^\prime$ of $G_{V^{\sh}}$ is $G^\prime_{V^{\sh}}$. We also point out that the $V$-schemes $G_1$, ..., $G_m=G^\prime$ are  of finite type.

\medskip\smallskip\noindent
{\bf 2.4. Representations.} 
We denote also by $\delta:R\to R\otimes_V R$ the comultiplication map of $G$. Let $L_0$ be a finite subset of $I_R$. For $l_0\in L_0$ we write $\delta(l_0)=\sum_{j\in I(l_0)} a_{0j}\otimes l_{0j}$, where $I(l_0)$ is a finite set and $a_{0j}$, $l_{0j}\in R$. Let $L$ be the $V$-submodule of $R$ generated by $1$, by $l_0$'s, and by $l_{0j}$'s. It is known that we have $\delta(L)\subset R\otimes_V L$, cf. [Ja, 2.13]. So $L$ is a $G$-module and thus we have a homomorphism $\rho(L):G\to GL(L)$ between flat, affine group schemes over $V$. Let $\scrB$ be a $V$-basis of $L$ contained in $\{1\}\cup I_R$. For $l\in \scrB\cap I_R$, we write $\delta(l)=\sum_{l^\prime\in\scrB} a_{ll^\prime}\otimes l^\prime$, where each $a_{ll^\prime}\in R$. As $\delta(l)-l\otimes 1+1\otimes l\in I_R\otimes_V I_R$, we have $a_{1l}=l$. But if $GL(L)=\Spec(A(L))$ and if $q(L):A(L)\to R$ is the $V$-homomorphism that defines $\rho(L)$, then $R(L):=\im(q(L))$ is the $V$-algebra generated by the $a_{ll^\prime}$'s. As a conclusion we have:

\medskip
{\bf (i)} The $V$-subalgebra of $R$ generated by $L$ is contained in $R(L)$. So if $G$ (resp. if $G_K$) is of finite type over $V$ (resp. over $K$) and if $L_0$ generates the $V$-algebra $R$ (resp. the $K$-algebra $R[{1\over\pi}]$), then $\rho(L)$ (resp. $\rho(L)_K$) is a closed embedding homomorphism. 

\medskip\smallskip
\noindent
{\bf 2.5. Quasi-sections.} Let $X$ be a reduced, flat $V$-scheme of finite type. Let $y\in X(k)$. From [Gr, Cor. (17.16.2)] we get the existence of a finite field extension $\tilde K$ of $K$ such that there is a faithfully flat, local $V$-subalgebra of $\tilde K$ that is of finite type and we have a morphism $z:\Spec(\tilde V)\to X$ whose image contains $y$. 

\medskip\noindent
{\bf 2.5.1. Lemma.} {\it {\bf (a)} If $V$ is complete, then we can assume $\tilde V$ is a discrete valuation ring.

\medskip
{\bf (b)} Let $a:Y\to X$ be a morphism between reduced, flat $V$-schemes of finite type. Suppose that $V$ is complete, that $k=\bar k$, and that for any $W\in\scrW$ the map $a(W):Y(W)\to X(W)$ is onto. Then the map $a(k):Y(k)\to X(k)$ is surjective.}

\medskip
\proof
As $V$ is complete, it is also a Nagata ring (cf. [Ma, (31.C), Cor. 2]) and thus the normalization $V^{\n}_{\tilde K}$ of $V$ in $\tilde K$ is a finite $V$-algebra. So as $\tilde V$ we can take any local ring of the normalization of $\tilde V$ in $\tilde K$ (it is a local ring of $V^{\n}_{\tilde K}$) that is a discrete valuation ring. From this (a) follows. By taking $W$ to be $\tilde V$, we get that $y\in\im(a(k))$. So as $y$ was arbitrary, we get that (b) holds.\endproof

\medskip\smallskip\noindent
{\bf 2.6. Lemma.} {\it Suppose $V$ is complete and $k=\bar k$. Let $f:F\to G$ be a homomorphism such that $f_K$ is an isomorphism. Then $f(V):F(V)\to G(V)$ is an isomorphism. If moreover $G$ is smooth, then $f$ is an isomorphism.}

\medskip
\proof      
Let $T$ be a maximal split torus of $F$, cf. [DG, Vol. III, Exp. XIX, 6.1]. We show that the assumption that $F(V)\lneqq G(V)$ leads to a contradiction. We have $F(K)=F(V)T(K)F(V)$, cf. Cartan decomposition of [BT, 4.4.3]. So as $F(V)\lneqq G(V)$, there is $g\in G(V)\cap (T(K)\setminus T(V))$. We write $T=\dbG_m^s=\Spec(V[w_1,...,w_s][{1\over{w_1...w_s}}])$. Let $w\in P$ be such that under the $V$-homomorphisms $P\twoheadrightarrow V[w_1,...,w_s][{1\over{w_1...w_s}}]\to K$ that define $g\in T(K)\leqslant F(K)$, it is mapped into an element of the set $\{w_i,w_i^{-1}|i\in\{1,...,s\}\}$ that maps into $K\setminus V$. As $f_K$ is an isomorphism we can identify $P[{1\over\pi}]=R[{1\over\pi}]$. Let $n\in\dbN$ be such that $\pi^nw\in R$. Under the $V$-homomorphism $R\to K$ that defines $g^{n+1}$, the element $\pi^nw$ maps into an element of $K\setminus V$. Thus $g^{n+1}\notin G(V)$. Contradiction. So $F(V)=G(V)$.  

Let now $G$ be smooth. We show that the assumption that $f$ is not an isomorphism leads to a contradiction. We can assume $k=\bar k$. As $R\neq P$, there is $w_0\in P\setminus R$ such that $x:=\pi$$w_0\in R\setminus \pi$$R$. Let $\bar g:R\to k$ be a $k$-homomorphism such that  $\bar g(x)\neq 0$. Let $g:R\to V$ be a $V$-homomorphism that lifts $\bar g$ (as $V$ is henselian and $G$ is smooth). The $K$-homomorphism $g[{1\over\pi}]:R[{1\over\pi}]\to K$ maps $w_0$ into an element of $K\setminus V$. So $g$ defines an element of $G(V)\setminus F(V)$. So $F(V)\lneqq G(V)$. Contradiction. Thus $f$ is an isomorphism.\endproof

\bigskip
\noindent
{\boldsectionfont \S3. Proof of Theorem 1.1}
\bigskip

In this chapter we prove 1.1. To prove 1.1 we can assume that $V$ is complete, that $k=\bar k$, and that $f_K$ is an isomorphism. So $f:F\to G$ is defined by a $V$-monomorphism $P\hookrightarrow R$ that induces a $K$-isomorphism $P[{1\over\pi}]\arrowsim R[{1\over\pi}]$ to be viewed as an identity. Let $\rho(L):G\to GL(L)$ be as in 2.4, with $L_0\in I_R\subset R$ such that it generates the $K$-algebra $P[{1\over\pi}]=R[{1\over\pi}]$. The generic fibre of $\rho(L)$ is a closed embedding, cf. 2.4 (i) and the choice of $L_0$. To prove 1.1 for $f$, it suffices to prove 1.1 for $\rho(L)\circ f$. So we can assume $\rho(L)$ is a closed embedding; so $G$ is a reduced, flat, closed subgroup of $GL(L)$ and so of finite type over $V$.

\medskip\smallskip\noindent
{\bf 3.1. Proofs of 1.1 (a) and (b).} 
Let $G^\prime=G_m\to ...\to G_1\to G_0=G$ be as in 2.3.2. As $F$ is smooth over $V$ and due to  2.3.2 (ii), $f:F\to G$ factors through a homomorphism $f^\prime:F\to G^\prime$. As $G^\prime$ is smooth, $f^\prime$ is an isomorphism (cf. 2.6). So 1.1 (a) holds. 

For $i\in\{0,...,m-1\}$, $G_i$ is a reduced, flat group scheme of finite type (cf. end of 2.3.2) and so a closed subgroup of some general linear group $H_i$ (cf. 2.4 (i)). So each group $\Ker(G_{i+1k}\to G_{ik})$ is a subgroup of a product of a finite number of copies of $\dbG_a$, cf. 2.3.1. As $f^\prime$ is an isomorphism, $\Ker(f_k)=\Ker(F_k\to G_k)$ has a composition series whose factors are subgroups of $\Ker(G_{i+1k}\to G_{ik})$ ($i\in\{0,...,m-1\}$). Thus $\Ker(f_k)$ is a unipotent group in the sense of [DG, Vol. II, Exp. XVII, 1.1]. As $\Ker(f_k)\vartriangleleft F_k$ and as $F_k$ has a trivial unipotent radical (being reductive), $\Ker(f_k)$ has dimension $0$. But $F_{k}$ is connected and so its action on $(\Ker(f_k))_{\red}$ via inner conjugation is trivial. So $(\Ker(f_k))_{\red}\leqslant Z(F_k)$. 

Let $\bar g\in (\Ker(f_k))(k)$. By induction on $i\in\{0,...,m\}$ we show that $\bar g\in\Ker(F(k)\to G_i(k))$. The case $i=0$ is obvious as $G_0=G$. For $i\in\{0,...,m-1\}$ the passage from $i$ to $i+1$ goes as follows. Let $S_i$ be the reduced subgroup of $G_{ik}$ such that $G_{i+1}$ is the dilatation of $G_i$ centered on $S_i$. Let $\tilde H_{i+1}$ be the dilatation of $H_i$ centered on $S_i$; we have a closed embedding homomorphism $G_{i+1}\hookrightarrow \tilde H_{i+1}$ (cf. [BLR, Prop. 2 (c) and (d), p. 64]). We consider a closed embedding homomorphism $\tilde H_{i+1}\hookrightarrow GL_{n_i}$ (with $n_i\in\dbN$, cf. 2.4 (i)). We have homomorphisms $F\to G_{i+1}\hookrightarrow \tilde H_{i+1}\hookrightarrow GL_{n_i}$. Let $\bar h\in GL_{n_i}(k)$ be the image of $\bar g$ via $F\to GL_{n_i}$. As $\bar g\in Z(F_k)(k)$, $\bar h$ is a semisimple element. As $\Ker(H_{i+1k}\to H_{ik})$ is a product of $\dbG_a$ groups (cf. 2.3.1) and as $\bar g\in\Ker(F(k)\to H_i(k))$, $\bar h$ is a unipotent element. So $\bar h$ is the identity element. So $\bar g\in\Ker(F(k)\to G_{i+1}(k))$. This ends the induction.

As $f^\prime:F\to G^\prime=G_m$ is an isomorphism, we get that $\bar g$ is the identity element of $F(k)$. Thus $(\Ker(f_k))_{\red}$ is a trivial group. So $\Ker(f_k)$ is also connected. So 1.1 (b) holds.

\medskip\smallskip\noindent
{\bf 3.2. Proof of 1.1 (c).} As $V$ is a Nagata ring, $R$ is also a Nagata ring (cf. [Ma, (31.H)]). So the normalization $G^{\n}=\Spec(R^{\n})$ of $G$ is a finite $G$-scheme. The homomorphism $F\to G$ factors through a morphism $F\to G^{\n}$. We have $F(V)=G(V)=G^{\n}(V)$, cf. 2.6; this also holds if $V$ is replaced by a $W\in\scrW$. Thus $f(k):F(k)\to G^{\n}(k)$ is an epimorphism, cf. 2.5.1 (b) applied to $F\to G^{\n}$. So $G^{\n}_k$ is connected and the morphism $F_k\to (G^{\n}_k)_{\red}$ is dominant. So both $F$ and $G^{\n}$ have unique local rings $O_1$ and $O_2$ (respectively) that are faithfully flat $V$-algebras and discrete valuation rings. The natural $V$-homomorphism $O_2\to O_1$ is dominant and becomes an isomorphism after inverting $\pi$. So we can identify $O_1=O_2$. So as $F_K=G^{\n}_K$, the normal, noetherian, affine schemes $F$ and $G^{\n}$ have the same set of local rings that are discrete valuation rings. Thus we have $R^{\n}=P$ (cf. [Ma, (17.H), Thm. 38]) and so $F\to G^{\n}$ is an isomorphism. So $f$ is a finite morphism. So 1.1 (c) holds.

\medskip\smallskip\noindent
{\bf 3.3. Proof of 1.1 (d).} As $\Ker(f_k)$ is unipotent, it is a subgroup of the unipotent radical of a Borel subgroup of some general linear group over $k$ (cf. [DG, Vol. II, Exp. XVII, 3.5]). So $\Lie(\Ker(f_k))$ is formed by nilpotent elements (cf. 2.2) and is normalized by $F_k$. The root data of $F_k$ and $F_{\bar K}$ are isomorphic (cf. [DG, Vol. III, Exp. XXII, 2.8]) and determine the isomorphism classes of $F_k$ and $F_{\bar K}$ (cf. [DG, Vol. III, Exp. XXIII, 5.1]). So the hypothesis of 1.1 (d) implies that either $p>2$ or $p=2$ and $F_k$ has no normal subgroup isomorphic to $SO_{2n+1}$ ($n\in\dbN$). So $\Lie(\Ker(f_k))$ has no non-trivial simple $F_k$-submodule, cf. 2.2.1. Thus $\Lie(\Ker(f_k))=0$. From this and the connectedness part of 1.1 (b), we get that $\Ker(f_k)$ is trivial. So $f_k$ is a closed embedding. As $f_k$ and $f_K$ are closed embeddings, from Nakayama's lemma we get that the finite morphism $F\times_{G} \Spec(\scrO)\to\Spec(\scrO)$ is a closed embedding for any local ring $\scrO$ of $G$. Thus $f$ is a closed embedding and so an isomorphism (as $f_K$ is so). So 1.1 (d) holds. This ends the proof of 1.1.

\medskip\smallskip\noindent
{\bf 3.4. Remarks.} {\bf (a)} The reference in [Va2, proof of 4.1.2] to [Va1, 3.1.2.1 c)] does not always work for $p=2$ (cf. 1.1 (d)). However, in [Va2, proof of 4.1.2] one can always choose $\rho_{B(k)}$ and $L$ such that $\rho$ is a closed embedding (cf. 2.4 (i)).

\smallskip
{\bf (b)} We continue to assume that $f_K$ is an isomorphism and that $k=\bar k$. As $\Lie(Z(F_k))$ is formed by semisimple elements, the connected, unipotent group $Z(F_k)\cap\Ker(f_k)$ is trivial. Suppose now that $p=2$ and $f$ is not an isomorphism. So $f_k$ is not a closed embedding (see 3.3) and so $(f_k)_{\red}:F_k\to (G_k)_{\red}$ is not a closed embedding. We get that $F_k^{\ad}\to (G_k)_{\red}/f_k(Z(F_k))$ is an isogeny (cf. 1.1 (c)) that is a finite product $\prod_{j\in J} f_j:F_j\to G_j$ of isogenies $f_j$ with $F_j$ as a simple, adjoint group. Any $f_j$ is a purely inseparable isogeny (cf. 1.1 (b)) and there is $j_0\in J$ such that $f_{j_0}$ is not an isomorphism. We have $\Ker(f_k):=\prod_{j\in J} \Ker(f_j)$ and so $\Ker(f_{j_0})$ is a non-trivial, unipotent subgroup of $F_j$ whose Lie algebra is formed by nilpotent elements (see 3.3).  From  2.2.1 and its proof we get that $F_{j_0}$ is an $SO_{2n+1}$ group and that $\Lie(\Ker(f_{j_0}))$ contains the unique simple $F_{j_0}$-submodule $\grn_{j_0}$ of $\Lie(F_{j_0})$ of dimension $2n$. The quotient of $F_{j_0}$  by $\grn_{j_0}$ (see [Bo, \S17]) is an $Sp_{2n}$ group (cf. [Bo, 23.6, p. 261]). So $f_{j_0}$ factors through a purely inseparable isogeny $g_{j_0}:Sp_{2n}\to G_{j_0}$ whose kernel is a unipotent group. So $\Lie(\Ker(g_{j_0}))$ is formed by nilpotent elements and so it is trivial (cf. 2.2.1). Thus $g_{j_0}$ is an isomorphism. As a conclusion, the root data of $F_k$ and $(G_k)_{\red}$ are not isomorphic and so $F_k$ and $(G_k)_{\red}$ are not isomorphic.

\bigskip
\noindent
{\boldsectionfont \S4. Proof of Theorem 1.2}
\bigskip
In this chapter we prove 1.2. To prove 1.2 we can assume that $V$ is complete, that $k=\bar k$, and that $\text{tr}.\text{deg}.(k)\Ge 1$.

\medskip\smallskip\noindent
{\bf 4.1. The group $\Gamma$.} As $k=\bar k$,  $(G_k)^0_{\red}$ is a split reductive group. Let $T_1$, ..., $T_s$ be a finite number of $\dbG_m$ subgroups of $(G_k)^0_{\red}$ that generate $(G_k)^0_{\red}$. For $i\in\{1,...,s\}$ let $y_i\in T_i(k)$ be an element of infinite order (here is the place where we need, in the case when $p\in\dbN$, that $\text{tr}.\text{deg}.(k)\Ge 1$). So the Zariski closure in $T_i$ of the subgroup of $T_i(k)$ generated by $y_i$, is $T_i$ itself. Let $\Gamma$ be the subgroup of $G(k)$ generated by $y_1$, ..., $y_s$. We conclude:

\medskip
{\bf (i)} The Zariski closure of $\Gamma$ in $G_k$ is $(G_k)^0_{\red}$.

\medskip\smallskip\noindent
{\bf 4.2. The finite type case.}
In this section we prove 1.2 under the assumption that $G$ is of finite type over $V$. For $i\in\{1,...,s\}$ let $V_i$ be a finite $V$-algebra that is a discrete valuation ring and such that there is $w_i\in G(V_i)$ that lifts $y_i$ (cf. 2.5.1 (a)). By replacing $V$ with its normalization in the composite field of the fields $V_i[{1\over\pi}]$ ($i\in\{1,...,s\}$), we can assume that for $i\in\{1,...,s\}$ there is $w_i\in G(V)$ that lifts $y_i$. Let $G^\prime$ be as in 2.3.2. As $w_1, ..., w_s\in G^\prime(V)=G(V)$, from 4.1 (i) we get that the Zariski closure of $\im(G^\prime_k(k)\to G_k(k))$ in $(G_k)_{\red}$ contains $(G_k)^0_{\red}$. So $(G_k)^0_{\red}\leqslant\im(G^\prime_k\to G_k)$. So if $G^{\prime0}_k$ is the identity component of $G^\prime_k$, then we have an isogeny $G^{\prime0}_k\to (G_k)^0_{\red}$. From this and [Bo, 14.11] we get that the unipotent radical of $G^{\prime0}_k$ is trivial. So $G^{\prime0}_k$ is a reductive group, cf. [Bo, 11.21].  

Let $G^{\prime0}$ be the open subgroup of $G^\prime$ formed by  $G^{\prime0}_k$ and by $G_K$. Let $(G_j^\prime)_{j\in J}$ be a covering of $G^\prime$ by open affine subschemes such that each $G^\prime_{jk}$ is an open subscheme of either $G^{\prime0}_k$ or $G^\prime_k\setminus G^{\prime0}_k$. The product scheme $G^{\prime0}\times_{G^\prime} G_j^\prime$ is: (i) $G_j^\prime$ if $G^\prime_{jk}\hookrightarrow G_k^{\prime0}$, and (ii) $G_{jK}^\prime$ if $G^\prime_{jk}\hookrightarrow G^\prime_k\setminus G_k^{\prime0}$. So $G^{\prime0}\times_{G^\prime} G_j^\prime$ is affine for any $j\in J$. Thus the morphism $G^{\prime0}\to G^\prime$ is affine and so $G^{\prime0}$ is affine. So $G^{\prime0}$ is a smooth, affine group scheme over $V$ whose special fibre is a reductive group and whose generic fibre is connected. This implies that $G^{\prime0}$ is a reductive group scheme over $V$, cf. [DG, Vol. III, Exp. XIX, 2.6 and 2.7]. So $G_K$ is a reductive group. So 1.2 (a) holds. 

Let $y\in G(k)$.  From 2.6 applied to $G^{\prime0}\to G^\prime$, we get that $G^{\prime0}=G^\prime$ and that $G^{\prime0}(V)=G^\prime(V)=G(V)$. So by replacing $V$ with a $W\in\scrW$, we can assume there is $w\in G(V)$ that lifts $y$ (cf. 2.5.1 (a)). So $y\in \im(G^{\prime0}(k)\to G(k))$. So we have an epimorphism $G^{\prime0}_k\twoheadrightarrow (G_k)_{\red}$ and so $G_k$ is connected. So 1.2 (a) holds.  We check 1.2 (b). From 1.1 (d) we get that the homomorphism $G^{\prime0}\to G$ is a closed embedding and so an isomorphism (as its generic fibre is so). So 1.2 (b) holds. We check that 1.2 (c) holds. 
The homomorphism $G^{\prime0}\to G^{\n}$ is an isomorphism (cf. 1.1 (c)) and so $G^{\n}$ is a reductive group scheme. So 1.2 (c) also holds. This ends the proof of 1.2 for the case when $G$ is of finite type.

\medskip\smallskip\noindent
{\bf 4.3. The general case.} To end the proof of 1.2 we are left to show that $G$ is of finite type. Let $I_0$ be the ideal of $R$ that defines $(G_k)^0_{\red}$. Let $\scrL$ be the set of all $G$--modules $L$ obtained as in 2.4. For $L\in\scrL$, let $L_0$, $L$, $R(L)$, and $\rho(L):G\to GL(L)$ be as in 2.4. We can identify $G(L):=\Spec(R(L))$ with the Zariski closure in $GL(L)$ of the image of $\rho(L)$. As $G_K$ is smooth over $K$ and connected, it is also of finite type. We now choose $L_0$ such that $\rho(L)_K$ is an isomorphism (cf. 2.4 (i)) and $L_0$ modulo $I_0$ generates the finite type $k$-algebra $R/I_0$. So $R(L)$ surjects onto $R/I_0$ (cf. 2.4 (i)) and thus $\rho(L)_k$ induces a closed embedding homomorphism $(G_k)^0_{\red}\hookrightarrow (G(L)_k)_{\red}^0$ between smooth, connected groups of dimension $\dim(G_K)$. So by reasons of dimensions, we get that $({G(L)_k})_{\red}^0$ is a reductive group isomorphic to $(G_k)^0_{\red}$. As in 4.2, by replacing $V$ with some $W\in\scrW$ we can assume that for any $i\in\{1,...,s\}$ there is $w_{iL}\in G(L)(V)$ that liftes $y_i\in ({G(L)_k})_{\red}^0(k)=(G_k)^0_{\red}(k)$. So from 4.2 applied to $G(L)$ (instead of $G$), we get that the normalization $G(L)^{\n}=\Spec(R(L)^{\n})$ of $G(L)$ is a reductive group scheme over $V$.

Let $\tilde L\in\scrL$ be such that $L\subset \tilde L$. So we have a homomorphism $\rho(\tilde L,L):G(\tilde L)\to G(L)$ such that $\rho(L)=\rho(\tilde L,L)\circ\rho(L)$. So $R(L)\hookrightarrow R(\tilde L)\hookrightarrow R\hookrightarrow R(L)[{1\over\pi}]=R(\tilde L)[{1\over\pi}]=R[{1\over\pi}]$. Let $\rho(\tilde L,L)^{\n}:G(\tilde L)^{\n}\to G(L)^{\n}$ be the morphism defined by $\rho(\tilde L,L)$. Let $W\in\scrW$ be such that there is $w_{i\tilde L}\in G(\tilde L)(W)$ that lifts $y_i\in {(G(\tilde L)_k})_{\red}^0(k)=(G_k)^0_{\red}(k)$. The normalization $(G(\tilde L)_W)^{\n}$ of $G(\tilde L)_W$ is a reductive group scheme over $W$, cf. 4.2 applied to $G(\tilde L)_W$. The morphisms $(G(\tilde L)_W)^{\n}\to G(\tilde L)_W^{\n}\to G(L)_W^{\n}\to G(L)_W$ define a homomorphism $\rho(\tilde L,L,W)^{\n}:(G(\tilde L)_W)^{\n}\to G(L)_W^{\n}$ between reductive group schemes over $W$ whose generic fibre is an isomorphism. So $\rho(\tilde L,L,W)^{\n}$ is an isomorphism, cf. 2.6. This implies that $\rho(\tilde L,L)^{\n}$ is an isomorphism. So $R(L)\hookrightarrow R(\tilde L)\hookrightarrow R(L)^{\n}=R(\tilde L)^{\n}\hookrightarrow R[{1\over\pi}]$. So as $R=\cup_{\tilde L\in\scrL} R(\tilde L)$, we have $V$-monomorphisms  $R(L)\hookrightarrow R\hookrightarrow R(L)^{\n}$. But $R(L^{\n})$ is a finite $R(L)$-algebra, cf. 4.2. So as $R(L)$ is a noetherian $V$-algebra, we get that $R$ is a finite $R(L)$-algebra and so a finitely generated $V$-algebra. This ends the proof of 1.2.

\bigskip
\noindent
{\boldsectionfont \S5. An application}
\bigskip 

We assume that $p\in\dbN$ is a prime and that $k=\bar k$. We take $V$ to be the Witt ring $W(k)$ of $k$. The goal of this section is to exemplify how one can use 1.1 (d) to extend results on semisimple groups over a field of characteristic $0$ to results on semisimple groups over $k$. Let $H_k$ be an absolutely simple, adjoint group over $k$. Let $Q_k$ be a parabolic subgroup of $H_k$ different from $H_k$. Let $U_k$ be the unipotent radical of $Q_k$ and let $L_k$ be a Levi subgroup of $Q_k$. 

\medskip\smallskip\noindent
{\bf 5.1. Proposition.}  {\it Let $\rho_k:L_k\to GL(\Lie(U_k))$ be the representation of the inner conjugation action of $L_k$ on $\Lie(U_k)$. Then $\rho_k$ is a closed embedding.}

\medskip
\proof
Let $H$ be the adjoint group scheme over $V=W(k)$ that lifts $H_k$. Let $Q$ be a parabolic subgroup of $H$ that lifts $Q_k$. Let $U$ be the unipotent radical of $Q$ and let $L$ be a Levi subgroup of $Q$ that lifts $L_k$. Let $\rho:L\to GL(\Lie(U))$ be the representation of the inner conjugation action of $L$ on $\Lie(U)$. Let $T$ be a maximal torus of $L$. Let $B$ be a Borel subgroup of $G$ such that $T\leqslant B\leqslant Q$. Let $\Lie(H)=\Lie(T)\oplus_{\alpha\in\Phi} \grg_{\alpha}$ be the Weyl decomposition of $\Lie(H)$ with respect to $T$. Let $\Delta=\{\alpha_1,...,\alpha_r\}$ be the basis of the root system $\Phi$ that corresponds to $B$; here $r\in\dbN$ is the rank of $H_k$. Let $\Phi_U$ be the subset of $\Phi$ such that $\Lie(U)=\oplus_{\alpha\in\Phi_U} \grg_{\alpha}$. For each $i\in\{1,...,r\}$, there is $\alpha\in\Phi_U$ that is the sum of $\alpha_i$ with an element of $\Phi_U\cup\{0\}$ which is a linear combination with coefficients in $\dbN\cup\{0\}$ of elements of $\Delta\setminus\{\alpha_i\}$. As $H$ is adjoint, this implies that the inner conjugation action of $T$ on $\Lie(U)$ is via characters of $T$ that generate the group of characters of $T$. Thus the restriction of $\rho$ to $T$ is a closed embedding. So the identity component of $\Ker(\rho_K)$ is a semisimple group over $K$ that has rank $0$. Thus $\Ker(\rho_K)$ is a finite, \'etale subgroup of $Z(H_K)$. As $Z(H_K)\leqslant T_K$ and as the intersection $T_K\cap \Ker(\rho_K)$ is trivial, we get that $\Ker(\rho_K)$ is trivial. So $\rho_K$ is a closed embedding. From 1.1 (d) we get that $\rho$ is a closed embedding, except perhaps when $p=2$ and $L_k$ has a normal subgroup that is an $SO_{2n+1}$ group for some $n\in\dbN$. So for the rest of the proof we can assume that $p=2$ and that $L_k$ has a normal subgroup $S_k$ that is an $SO_{2n+1}$ group for some $n\in\dbN$, $n\Le r$. 

This implies that $r\Ge 2$ and that $H_k$ has a subgroup normalized by $T_k$ and which is a $PGL_2=SO_3$ group. So $H_k$ is an $SO_{2r+1}$ group, cf. [Va2, 3.8]. If $\Lie(\Ker(\rho_k))=\{0\}$, then as in the end of 3.3 we argue that $\rho$ is a closed embedding. So to end the proof, we only need to show that the assumption that $\Lie(\Ker(\rho_k))\neq \{0\}$ leads to a contradiction. As in 3.3 we argue that 1.1 (a) implies that $\Lie(\Ker(\rho_k))$ is formed by nilpotent elements. Based on 2.2.1 and its proof, we can assume that $S_k$ is such that $\Lie(S_k)\cap \Lie(\Ker(\rho_k))$ contains the ideal $\grn$ of $\Lie(S_k)$ generated by eigenspaces of the adjoint action of $T$ on $\Lie(H^{\ad})$ that correspond to short roots; we have $\dim_k(\grn)=2n$ (cf. proof of 2.2.1). Let $U_k^-$ be the unipotent subgroup of $H_k$ that is the opposite of $U_k$ with respect to $T_k$; so $\Lie(U_k^-)=\oplus_{\alpha\in\Phi_U} \grg_{-\alpha}\otimes_V k$. Let $w_0\in H_k(k)$ be such that it normalizes both $T_k$ and $L_k$ and we have $w_0U_kw_0^{-1}=U_k^-$, cf. [Bou, (XI), PLATE II].  As $\grn\subset\Lie(\Ker(\rho_k))$, $\grn$ centralizes $\Lie(U_k)$. As $\grn$ is a characteristic ideal of $S_k$, $\grn$ is normalized by $w_0$ and so $\grn$ also centralizes $\Lie(U_k^-)$. Moreover, $\grn$ is normalized by $\Lie(L_k)$. Thus $\grn$ is normalized by $\Lie(H_k)=\Lie(U_k)\oplus\Lie(L_k)\oplus\Lie(U_k^-)$ and so it is an ideal of $\Lie(H_k)$. But $\Lie(H_k)$ has a unique minimal ideal that has dimension $2r$, cf. [Hu, ($B_r$) of 0.13]. Thus $2n=\dim_k(\grn)\Ge 2r\Ge 2n$. Thus $n=r$ and so $S_k=H_k$. This implies that $Q_k=H_k$ and so we reached a contradiction to the relation $Q_k\neq H_k$.\endproof
\bigskip
\references{37}
{\nspace{

\Ref[BLR]
S. Bosch, W. L\"utkebohmert, M. Raynaud,
\sl N\'eron models,
\rm Springer-Verlag, 1990.

\Ref[Bo]
A. Borel,
\sl Linear algebraic groups,
\rm Grad. Texts in Math., Vol. {\bf 126}, Springer-Verlag, 1991.

\Ref[Bou]
N. Bourbaki,
\sl Lie groups and Lie algebras, 
\rm Chapters {\bf 4--6}, Springer-Verlag, 2002.

\Ref[BT]
F. Bruhat and J. Tits, 
\sl Groupes r\'eductifs sur un corps local: I. Donn\'ees radicielles valu\'ees,
\rm Inst. Hautes \'Etudes Sci. Publ. Math., Vol. {\bf 41}, pp. 5--251, 1972.

\Ref[DG]
M. Demazure, A. Grothendieck, \'et al.,
\sl Sch\'emas en groupes. Vol. I-III
\rm Lecture Notes in Math., Vol. {\bf 151--153}, Springer-Verlag, 1970.

\Ref[Gr] 
A. Grothendieck, 
\sl \'El\'ements de g\'eom\'etrie alg\'ebrique. IV. \'Etude locale des sch\'emas et des morphismes de sch\'ema (Quatri\`eme Partie), 
\rm Inst. Hautes \'Etudes Sci. Publ. Math., Vol. {\bf 32}, 1967.

\Ref[Hi]
G. Hiss,
\sl Die adjungierten Darstellungen der Chevalley-Gruppen, 
\rm Arch. Math. {\bf 42} (1982), pp. 408--416. 

\Ref[Hu]
J. E. Humphreys, 
\sl Conjugacy classes in semisimple algebraic groups, 
\rm Math. Surv. and Monog., Vol. {\bf 43}, Amer. Math. Soc., 1995.

\Ref[Ja]
J. C. Jantzen,
\sl Representations of algebraic groups. Second edition,
\rm Math. Surveys and Monog., Vol. {\bf 107}, Amer. Math. Soc., 2000.

\Ref[Ma] H. Matsumura, 
\sl Commutative algebra. Second edition, 
\rm The Benjamin/Cummings Publ. Co., Inc., Reading, Massachusetts, 1980.

 \Ref[Pi]
R. Pink,
\sl Compact subgroups of linear algebraic groups,
\rm J. of Algebra {\bf 206} (1998), pp. 438--504.

\Ref[PY]
G. Prasad and J.-K. Yu, 
\sl On quasi-reductive group schemes,
\rm math.NT/0405381, 34 pages revision, June 2004.

\Ref[Va1]
A. Vasiu,
\sl Integral canonical models of Shimura varieties of preabelian type,
\rm Asian J. Math., Vol. {\bf 3} (1999), no. 2, pp. 401--518.

\Ref[Va2] 
A. Vasiu,
\sl Surjectivity criteria for $p$-adic representations, Part I,
\rm Manuscripta Math. {\bf 112} (2003), no. 3, pp. 325--355.

}}

\noindent
$\vfootnote{} {Adrian Vasiu, Mathematics Department, University of Arizona, 617 N. Santa Rita, P.O. Box 210089, Tucson, AZ-85721, USA. e-mail: adrian\@math.arizona.edu}$
\enddocument